\documentclass{amsart}
\newtheorem{theorem}{Theorem}[section]

\newtheorem{lemma}[theorem]{Lemma}
\newtheorem{corollary}[theorem]{Corollary}

\newtheorem{remark}[theorem]{Remark}
\newtheorem{conjecture}{Conjecture}
\newtheorem{proposition}[theorem]{Proposition}

\newtheorem{definition}[theorem]{Definition}

\newcounter{rmnum}

\numberwithin{equation}{section}

\begin{document}
\title{Proper  actions on   topological groups: applications to quotient spaces}
\author{Sergey A. Antonyan }

\address{Departamento de  Matem\'aticas,
Facultad de Ciencias, Universidad Nacional Aut\'onoma de M\'exico,
 04510 M\'exico Distrito Federal,  M\'exico.}
\email{antonyan@unam.mx}

\begin{abstract}

Let $X$ be a Hausdorff topological group and $G$ a locally compact
subgroup of $X$. We show  that the natural action of $G$ on $X$ is proper in the sense of R.~Palais.
This is applied to prove that  there exists a closed set
$F\subset X$ such that $FG=X$ and the restriction   of the quotient projection  $X\to X/G$ to $F$ is a perfect map $F\to X/G$. This is a key
result to  prove  that   many topological properties (among them,
paracompactness and normality)  are transferred   from $X$ to
$X/G$, and some others are transferred  from $X/G$ to  $X$. Yet
another application leads to the inequality   $ {\rm dim}\, X\le
{\rm dim}\, X/G + {\rm dim}\, G$ for every paracompact topological group $X$
and a  locally compact subgroup $G$ of $X$  having  a compact group of connected components.
\end{abstract}

\thanks {The  author was supported in part by grants \#IN102608 from PAPIIT (UNAM) and  \#79536 from CONACYT (Mexico)}

\maketitle
\markboth{SERGEY A. ANTONYAN}{APPLICATIONS OF PROPER ACTIONS}

\section{Introduction}

By a $G$-space we mean a completely regular Hausdorff space together with a fixed continuous action of a Hausdorff topological  group  $G$ on it.

The notion of a proper $G$-space  was introduced in 1961 by R.~Palais \cite{pal:61} with the purpose to extend a substantial portion of the   theory of compact Lie group actions to the case of noncompact ones.

A $G$-space   $X$ is called proper (in the sense of Palais \cite[Definition 1.2.2]{pal:61}), if each point of $X$ has a, so called, {\it small} neighborhood, i.e., a neighborhood $V$ such that for every point of $X$ there is a neighborhood $U$ with the property that the set
$\langle U,V\rangle=\{g\in G \ | \  gU\cap V\not= \emptyset\}$    has compact closure in $G$.

 Clearly, if $G$ is compact, then  every $G$-space is  proper.

Many important problems in the theory of proper actions are
conjugated (see \cite{haj:71}, \cite{ab:74}, \cite{ab:78},
\cite{anne:03}, \cite{ant:05}) to the following major open
problem:

\begin{conjecture}\label{Conj1} Let $G$ be a locally compact group. Then the orbit space $X/G$ of any paracompact proper $G$-space $X$ is paracompact.
\end{conjecture}

This conjecture is open even if $X$ is metrizable; in this case it is equivalent  (see \cite{anne:03}) to the following old  problem going back to  R. Palais \cite{pal:61}:

\begin{conjecture}\label{Conj2} Let $G$ be a locally compact group and $X$ a metrizable proper $G$-space. Then the topology of $X$ is metrizable  by a  $G$-invariant metric.
\end{conjecture}

Due to Palais~\cite{pal:61}, it is known that Conjecture~\ref{Conj2} is  true  for a separable metrizable proper $G$-space $X$  provided  the acting group  $G$ is  Lie. Other special  cases are discussed in \cite{koz:65} and \cite{anne:03}. In particular, in \cite{anne:03},  Palais' result is   extended to  the case of a locally compact separable $G$ and a  metrizable locally separable $X$.

\

In this  paper we prove  Conjecture~\ref{Conj1} in  an important special  case, namely, when  $X$ is a  topological group endowed with the  natural action of its locally compact subgroup $G$ (see Corollary~\ref{C:2}). We first show that    $X$ is a proper $G$-space and then we
 establish  a more  general result (Theorem~\ref{T:12}) which has many interesting applications  in the theory of topological groups.

Below all topological groups are assumed to satisfy the Hausdorff separation axiom.

\begin{theorem}\label{T:11}  Let $X$ be a topological  group and $G$ a locally compact subgroup of $X$. Then the action of $G$ on $X$ given by the formula $g*x=xg^{-1}$, $g\in G$,  $x\in X$, is proper.
\end{theorem}

Recall that a subset $S$ of a proper $G$-space $X$ is called {\it $G$-fundamental}, or just {\it fundamental}, if  $S$ is a small set  and the saturation $G(S)=\{gs \ | \ g\in G, s\in S\}$ coincides with $X$.

Here is the key result of the paper:

\begin{theorem}\label{T:12}  Let $X$ be a topological  group and $G$ a locally compact subgroup of $X$. Then there exists a closed  $G$-fundamental set in $X$.
\end{theorem}

It is easy to prove (see \cite[Proposition~1.4]{ab:78}) that in every proper $G$-space $X$, the restriction of the orbit map $p:X\to X/G$ to any closed small set is perfect (i.e., is closed and has compact fibers). In combination with Theorem~\ref{T:12} this yields the following:

\begin{corollary}\label{C:0}  Let $X$ be a topological  group,  $G$ a locally compact subgroup of $X$, and $X/G$  the quotient space of all right cosets  $xG=\{xg\ | \ g\in G\}$, $x\in X$.
Then there exists a closed subset $F\subset X$ such that the restriction  $p|_F :F\to X/G$ is a perfect surjective  map.
\end{corollary}

This fact has the following immediate corollary about transfer of properties from $X$ to $X/G$:

\begin{corollary}\label{C:1}  Let $\mathcal P$ be a topological property stable under perfect  maps and also inherited by closed subsets.
 Assume that $X$ is a topological group with the  property  $\mathcal P$ and let $G$ be a locally compact subgroup of $X$. Then the quotient space  $X/G$ also has the property $\mathcal P$.
\end{corollary}

Among  such properties $\mathcal P$ we single out just some of those  which provide new results in Corollary~\ref{C:1}; these are: paracompactness, countable paracompactness,  weak paracompactness,  normality, perfect normality, \v Cech-completeness,  being a $k$-space (see~\cite[\S 5.1, \S 5.2, \S 5.3, \S 1.5, \S 3.9,  \S 3.3]{eng:77}) and stratifiability (see \cite{bor:66}). Thus,  we get the following positive solution of Conjecture~1 in an important  special  case:

\begin{corollary}\label{C:2}  Let  $X$ be a paracompact  topological group   and $G$  a locally compact subgroup of $X$. Then the quotient space  $X/G$ is paracompact.
\end{corollary}

In this connection it is in order to recall the following remarkable result of A.~V.~Arhangel'skii~\cite{arh:81}: every topological group is the  quotient of a paracompact zero-dimensional group. Hence, the  local compactness of $G$ is essential in Corollary~\ref{C:2}.

\medskip

We recall that a locally compact group  is called {\it almost connected}  if  its space of connected components endowed with the quotient topology is compact.

\begin{corollary}\label{C:3}  Let $\mathcal P$ be a topological property stable under open perfect   maps  and also inherited by closed subsets.  Assume that $X$ is a paracompact  group with the  property  $\mathcal P$ and let $G$ be an almost connected  subgroup of $X$. Then the quotient space  $X/G$ also has the property $\mathcal P$.
\end{corollary}

Among properties stable under open perfect   maps  and also inherited by closed subsets we highlight  strong paracompactness and realcompactness (see \cite[Exercises~5.3.C(c),  5.3H(d), and Theorem~3.11.4 and Exercises~3.11.G] {eng:77}). Thus,  we get the following:

\begin{corollary}\label{C:4}  Let  $X$ be a strongly paracompact (resp., paracompact and realcompact) topological group   and $G$  an almost connected  subgroup of $X$. Then the quotient space  $X/G$ is strongly paracompact (resp., paracompact and realcompact).
\end{corollary}

\begin{remark}
We note that the converse of the first statement in this corollary is not true. Namely, it is known that the Baire space $ B(\aleph_1)$ of weight $\aleph_1$ (which is, in fact, homeomorphic to  a commutative metrizable topological group) is strongly paracompact and its product with the additive group $\Bbb R$ of the reals is not (see \cite[Exercises~5.3F(a) and 5.3F(b)]{eng:77}). Hence, the direct product  $X=\Bbb R\times B(\aleph_1)$ and its subgroup $G=\Bbb R\times\{0\}$ provide the desired counterexample, answering negatively a question  from \cite{arhusp}. Furthermore, \cite[Open Problem 3.2.1, p.151]{arhtk} asks whether every  locally strongly paracompact group is strongly paracompact? The same group $\Bbb R\times B(\aleph_1)$ provides  a negative answer to this question too. Indeed, since the product of a compact space and a strongly paracompact space  is strongly paracompact (see \cite[Exercise~5.3H(a)]{eng:77}), we infer that $\Bbb R\times B(\aleph_1)$ is locally strongly paracompact.
\end{remark}

\medskip

Combining our Corollary~\ref{C:2} with a result of Abels~\cite[Main Theorem]{ab:74}, we obtain the following:

\begin{corollary}  Let   $X$ be a paracompact  group, $G$ an almost connected  subgroup of $X$, and  $K$  a maximal compact subgroup of $G$.    Then there exists a $K$-invariant subset  $S\subset X$ such that $X$ is $K$-homeomorphic to the product $G/K\times S$. In particular, $X$ is homeomorphic to $\Bbb R^n\times S$ for some $n\ge 1$.
\end{corollary}

In \cite{arh:05}  A.~V.~Arhangel'skii has studied  properties which
are transferred in the opposite direction, i.e., from $X/G$ to
$X$. The next corollary is a unified result of this sort which
implies many of those in \cite{arh:05} as well as provides some
new ones.

\begin{corollary}\label{C:5}  Let $\mathcal P$ be a topological property invariant  and inverse invariant of perfect  maps,  and also stable under multiplication  by a locally compact group.
 Assume that $X$ is a topological   group and let $G$ be a locally compact subgroup of $X$ such that  the quotient space $X/G$ has  the   property  $\mathcal P$. Then   the group $X$  also has the property $\mathcal P$.
\end{corollary}

Among  such properties  we mention  just some:   paracompactness, being a $k$-space, \v Cech-completeness  (see~\cite[\S 5.1, \S 3.3, \S 3.9]{eng:77}). That paracompactness is stable under multiplication  by a locally compact group follows from a result of Morita~\cite{mor:53} since every  locally compact group is paracompact  (even, strongly paracompact~\cite[Theorem~3.1.1]{arhtk}).

\begin{corollary}\label{C:6}  Let $\mathcal P$ be a topological property invariant  and inverse invariant of open perfect  maps,  and also stable under multiplication  by a locally compact group.
 Assume that $X$ is a paracompact  group and let $G$ be an almost connected  subgroup of $X$ such that  the quotient space $X/G$ has  the   property  $\mathcal P$. Then   the group $X$  also has the property $\mathcal P$.
\end{corollary}

Among such properties  we highlight   realcompactness (see
\cite[Theorem~3.11.14 and  Exercise~3.11.G, and also take into
account  that every locally compact group is realcompact]
{eng:77}).

\smallskip

Coorollary\ref{C:2} is further applied to prove the following Hurewicz type formula:

\begin{theorem}\label{T:13} Let $X$ be a  paracompact  topological group and   $G$  an almost connected  subgroup of   $X$.
   Then
$$ {\rm dim}\, X\le {\rm dim}\, X/G + {\rm dim}\, G.$$
\end{theorem}

\begin{remark}[\cite{skl:64}]\label{R} If in this theorem $X$ is a locally compact group then, in fact, the equality holds:
$$ {\rm dim}\, X = {\rm dim}\, X/G + {\rm dim}\, G.$$
\end{remark}

All the proofs are given in section~3.
\medskip

\section{Preliminaries}

Throughout the paper, unless otherwise is stated, by  a {\it
group} we shall mean   a topological group $G$ satisfying the
Hausdorff separation axiom; by $e$ \ we shall denote the unity of
\ $G$.

All topological spaces  are assumed to be Tychonoff (= completely
regular and Hausdorff). The basic ideas and facts of the theory of
$G$-spaces or topological transformation groups can be found in
G.~Bredon \cite{br:72} and in  R.~Palais \cite{pal:60}. Our basic
reference on  proper group actions is Palais' article
\cite{pal:61}.  Other good sources are  \cite{koz:65}, \cite{ab:74} and   \cite{ab:78}.

For the  convenience of the reader we recall, however,  some more
special definitions and facts below.

\smallskip

By a $G$-space we mean a topological space $X$ together with a fixed continuous action  $G\times X\to X$ of a
topological  group $G$ on $X$. By $gx$ we shall denote the image of the pair $(g, x)\in G\times X$  under the
action.

 If $Y$ is another $G$-space, a
continuous map $f:X\to Y$ is called a $G$-map or an equivariant
map, if $f(gx)=gf(x)$ for every $x\in X$ and $g\in G$.

If $X$ is a $G$-space, then for a subset $S\subset X$ and for a subgroup $H\subset G$, \ the $H$-hull (or
$H$-saturation) of $S$ is defined as follows:  $H(S)$= $\{hs \ |\ h\in H,\ s\in S\}$. If $S$ is the  one point set
$\{x\}$, then the $G$-hull $G(\{x\})$ usually is denoted by $G(x)$  and called  the orbit of $x$. The orbit
space $X/G$ is always considered in its quotient topology.

A subset $S\subset X$ is called $H$-invariant if it coincides with
its $H$-hull, i.e., $S=H(S)$. By an invariant set we shall mean a
 $G$-invariant set.

For any $x\in X$, the subgroup   $G_x =\{g\in G \ |  \  gx=x\}$ is
called  the stabilizer (or stationary subgroup) at $x$.

 A compatible metric $\rho$ on a metrizable $G$-space $X$  is called invariant or $G$-invariant, if
$\rho(gx, gy)=\rho(x, y)$ for all $g\in G$ and $x, y\in X$. If $\rho$ is a $G$-invariant
metric on any $G$-space  $X$, then it is easy to verify that  the formula
$$\widetilde{\rho}\bigl(G(x), G(y)\bigr)=inf\{\rho(x', y')  \ | \  \ x'\in G(x), \ y'\in G(y)\}$$
 defines  a pseudometric $\widetilde{\rho}$, compatible with the quotient topology of $X/G$. If, in addition, $X$ is a proper $G$-space then  $\widetilde{\rho}$ \ is, in fact,  a compatible metric on $X/G$ \cite[Theorem~4.3.4]{pal:61}.

For a closed subgroup $H \subset G$, by $G/H$ we will denote the $G$-space of cosets $\{gH | \ g\in G\}$ under
the action induced by left translations.

A locally compact group $G$ is called {\it almost connected} if the quotient  group  $G/G_0$ of $G$ modulo the connected component $G_0$ of the identity is compact.

 Such a group has a maximal compact subgroup $K$, i.e., every compact subgroup of $G$ is
  conjugate to a subgroup of $K$ \cite[Theorem A.5]{ab:74}. The corresponding classical theorem on Lie groups can be found in
   \cite[Ch.~XV, Theorem~3.1]{hoch:65}.

\medskip
\normalfont

 In 1961 Palais~\cite{pal:61} introduced the very important concept of a {\it
proper action} of an arbitrary locally compact  group $G$  and extended   a
substantial part of the theory of compact Lie transformation
groups to noncompact  ones.

 Let  $X$ be  a $G$-space. Two subsets  $U$ and $V$ in  $X$  are called  thin relative to each other
  \cite[Definition 1.1.1]{pal:61},  if the set
$\langle U,V\rangle=\{g\in G | \ gU\cap V\ne \emptyset\}$, called {\it the transporter} from $U$ to $V$,
has a compact closure in $G$.
   A subset $U$ of a $G$-space $X$ is called {\it $G$-small}, or just {\it small}, \ if  every point in $X$ has a neighborhood thin relative to $U$. A $G$-space $X$
is called  {\it  proper} (in the sense of Palais),  if   every point in  $X$ has a small neighborhood.

Clearly,  if $G$ is compact, then  every $G$-space is  proper. Furthermore, if $G$ acts properly on
a compact space,  then $G$ has to be compact as well. If  $G$ is discrete and $X$ is locally
compact, the notion of a proper action is the same as the classical notion of a {\it properly
discontinuous} action. When $G$=$\Bbb R$, the additive group of the reals,  proper $G$-spaces are
precisely the  {\it dispersive} dynamical systems with regular orbit space (see \cite[Ch.
IV]{basz:70}).

Important examples of  proper $G$-spaces are the coset spaces
$G/H$ with $H$ a compact subgroup of a locally compact group $G$.
Other interesting examples the reader can find in \cite{ab:74},
\cite{ab:78}, \cite{ant:98}, \cite{ant:99} and \cite{koz:65}.

\smallskip

In what follows we shall need also  the definition of a twisted product $G\times_K S$, where
 $K$ is a  closed subgroup of $G$,  and   $S$  a $K$-space.
$G\times_KS$ is the orbit space of the $K$-space $G\times S$ on which  $K$ acts by the rule: $k(g,
s)=(gk^{-1}, ks)$. Furthermore, there is a natural action of $G$ on $G\times_K S$ given by
$g^\prime[g, s]=[g^\prime g, s]$, where $g'\in G$ and $[g, s]$ denotes the $K$-orbit of the point
$(g, s)$ in $G\times S$.  We shall identify $S$, by means of the $K$-equivariant embedding $s\mapsto [e, s]$, $s\in S$,  with  the $K$-invariant
subset $\{[e, s] \ | \ s\in S\}$ of $G\times_KS$. This  $K$-equivariant embedding $S\hookrightarrow G\times_KS$ induces a homeomorphism of the $K$-orbit space $S/K$ onto the $G$-orbit space $(G\times_KS)/G$ (see \cite[Ch.~II, Proposition~3.3]{br:72}).

 The twisted products  are  of a particular interest in the theory of
transformation groups  (see \cite[Ch.~II, \S~2]{br:72}). It turns
out that every $G$-space  locally is a twisted product. For a more
precise formulation we  need to  recall  the following well known
notion  of a slice (see \cite[p.~305]{pal:61}):

\begin{definition}\label{D:21} Let $X$ be a $G$-space and   $K$   a closed  subgroup of $G$.
 A  $K$-invariant  subset $S\subset X$ is called a $K$-kernel if  there is a $G$-equivariant map $f:G(S)\to G/K$, {called the slicing map,}
 such that $S$=$f^{-1}(eK)$.   The saturation $G(S)$ is called   a {\it tubular} set and the subgroup $K$ will be
 referred as
   the slicing subgroup.

If in addition   $G(S)$ is open in $X$ then  we shall call $S$ a
$K$-slice in $X$.

  If  $G(S)=X$ then  $S$ is called {\it a global} $K$-slice of $X$.
\end{definition}

It turns out that each tubular set with a compact slicing subgroup
is a twisted product. The tubular neighborhood $G(S)$ is
$G$-homeomorphic to the twisted product $G\times_K S$; namely the
map  $\xi:G\times_K S\to G(S)$ defined by $\xi([g, s])=gs$ is a
$G$-homeomorphism (see \cite[Ch. II, Theorem 4.2]{br:72}). In what
follows we will  use this fact without a specific reference.

 One of the most powerful results in  the theory of topological
transformation groups is Palais' slice theorem~\cite[Proposition 2.3.1]{pal:61} which states  that, if $X$ is a proper $G$-space with $G$ a Lie
 group, then for any point $x\in X$, there exists a $G_x$-slice $S$ in $X$ such that   $x\in S$.
 In general, when $G$ is not a Lie group, it is no longer true that a $G_x$-slice exists at each point of $X$ (see
\cite{ant:94}). Generalizing the case of Lie group actions, in  \cite{ab:78} and \cite{ant:05} (see also \cite{ant:94} for the case of compact non-Lie group actions),  approximate versions of Palais' slice theorem  for non-Lie group actions were proved. Below, in  the proof of
Theorem~\ref{T:13}, we shall need the following global slice theorem established by  H.~Abels~ \cite[Main Theorem]{ab:74}:

\begin{theorem}[Global Slice  Theorem] \label{T:GSlice} Let $G$ be an almost connected  group, $K$ a maximal compact subgroup of $G$, and $X$  a proper $G$-space with   a paracompact orbit space. Then $X$ admits a global $K$-slice.
 \end{theorem}

\medskip

On any group $G$ one  can define two natural  (but equivalent) actions of $G$ given  by the formulas
$$g\cdot x=gx,  \quad \quad\text{and}\quad g*x=xg^{-1},$$
respectively, where in the right parts the group operations are used  with $g, x\in G$.

  Throughout we shall consider the second action only.

\

  By $U(G)$ we shall denote the Banach space of all right uniformly continuous bounded functions $f:G\to \Bbb R$ endowed with the supremum norm. Recall that $f$ is called right  uniformly continuous, if for every $\varepsilon >0$ there exists a neighborhood $O$ of the unity in $G$ such that $|f(y)-f(x)|<\varepsilon$ whenever $yx^{-1}\in O$.

  We shall consider the induced action of $G$ on $U(G)$, i.e.,
  $$(gf)(x)=f(xg), \quad\text{for all}\quad  g, x\in G.$$
  It is easy to check that this action is continuous, linear  and isometric (see e.g., \cite[Proposition~7]{ant:87}).

\begin{proposition}\label{P1}  Let $G$ be a  group. Then for  every  $f\in U(G)$,  the map
$$f_*:G\to U(G)$$
defined by $ f_*(x)(g)=f(xg^{-1}),\  x, g\in G $, is a right uniformly continuous $G$-map.
\end{proposition}
\begin{proof} A simple verification.
\end{proof}

\medskip

\section{Proofs}

\noindent
{\it Proof of Theorem~\ref{T:11}}. Choose a neighborhood $U$ of the identity in $X$ such that $U=U^{-1}$ and $U^3\cap G$ has a compact closure in $G$. We claim that for every $x\in X$, the neighborhood $xU$ is $G$-small. Indeed, let $y\in X$ be any point. Two cases are possible.

{\it Case 1}. Assume that $y\in xU^2G$. Then  $y=xu_1u_2h$ with $u_1, u_2\in U$ and $h\in G$. We claim that  $xU^2h$ is a neighborhood of $y$ thin relative to $xU$. Indeed, if  $g\in \langle xU, xU^2h\rangle$, then  $g^{-1}h^{-1}\in U^3\cap G$. Since $U^3\cap G$ has a compact closure, we see that so does  the set  $(U^3\cap G)h$  which contains  $g^{-1}$. This yields that $\langle xU, xU^2h\rangle$ is contained in   $h^{-1}(U^3\cap G)^{-1}$, which also has a compact closure. Hence the transporter $\langle xU, xU^2h\rangle$ has a compact closure, as required.
\medskip

{\it Case 2}.  Assume that $y\notin xU^2G$. In this case $y\notin \overline{xUG}$. Indeed, if  $y\in \overline{xUG}$ then the neighborhood $xUx^{-1}y$ of  $y$ should meet the set $xUG$. Then $xux^{-1}y=xvh$ for suitable elements $u, v\in U$ and $h\in G$. Then $y=xu^{-1}vh\in xU^2G$, a contradiction.

Hence the open set $V=X\setminus  \overline{xUG}$ is a  $G$-invariant  neighborhood of $y$, and $V$  is thin relative to $xU$ because the transporter
$ \langle xU, V\rangle$ is empty in this case. \qed

\medskip

 \begin{lemma}\label{P:right}
Let $U$ be a unity  neighborhood in $X$ such that  $U=U^{-1}$ and $U^3\cap G$ has compact closure in $X$. Then for every $x\in X$, the neighborhood $Ux$ is a $G$-small set  (under the action $g*x$).
 \end{lemma}
\begin{proof}
Take $y\in X$ arbitrary. Since by Theorem~\ref{T:11}, $xU$ is a $G$-small set, the point $xy$ must have a neighborhood $V$ which is thin relative to $xU$.
We claim that the neighborhood $x^{-1}V$ of the point $y$ is the desired one, i.e., it is thin relative to $Ux$. This claim will follow from the following equality:
$$\langle Ux, x^{-1}V\rangle= \langle xU, V\rangle x.$$
In its turn,  this  equality results from  the following chain of obvious equivalences:
  \begin{align*}
  g\in \langle Ux, x^{-1}V\rangle &\Longleftrightarrow   Uxg^{-1}\cap x^{-1}V\ne\emptyset  \Longleftrightarrow  x\big(Uxg^{-1}\cap x^{-1}V\big)\ne\emptyset\\
    & \Longleftrightarrow xUxg^{-1}\cap V\ne\emptyset \Longleftrightarrow gx^{-1}\in\langle xU, V\rangle  \Longleftrightarrow g\in\langle xU, V\rangle x.
    \end{align*}

  \smallskip

  Thus, the closure $\overline{\langle Ux, x^{-1}V\rangle}$ equals to $\big(\overline{\langle xU, V\rangle}\big)x$. Since $\overline{\langle xU, V\rangle}$ is compact and the right translations are autohomeomorphisms of $X$, we conclude that $\big(\overline{\langle xU, V\rangle}\big)x$, and hence,  $\overline{\langle Ux, x^{-1}V\rangle}$ \   is compact.  This completes the proof.

\end{proof}

\begin{proposition}\label{P:31}  Let $X$ be a  group and $G$ a locally compact subgroup of $X$. Then there exists a locally finite covering of $X$ consisting of $G$-invariant open sets of the form $S_iG$, where each $S_i$ is an open  $G$-small subset of $X$.
\end{proposition}
\begin{proof}  By Theorem~\ref{T:11}, $X$ is a proper $G$-space, and hence, one can choose a  $G$-small neighborhood $U$ of the unity in  $X$.  By virtue of Markov's theorem~\cite[Theorem~3.3.9]{arhtk}, there exists a   right  uniformly continuous function  $f:X\to [0, 1]$  such that
 \begin{equation}\label{1}
 f(e)=0 \quad \text{and}\quad f^{-1}\big([0, 1)\big)\subset U.
\end{equation}

Then, by Proposition~\ref{P1}, $f$ induces an  $X$-equivariant map
$f_*:X\to U(X)$ defined by the rule:
$$ f_*(x)(g)=f(xg^{-1}),\quad x, g\in X. $$

Denote by $Z$ the image $f_*(X)$. Clearly, $Z$ is the $X$-orbit of the
point $f_*(e)$ in the $X$-space $U(X)$, and  the metric of $U(X)$ induces an  $X$-invariant
metric on $Z$. We claim that
\begin{equation}\label{2}
f_*^{-1}(\Gamma_{x, V})\subset x^{-1}*U=Ux, \quad\text{for every}\quad
x\in X,
\end{equation}
where $V=[0, 1)$ and $\Gamma_{x, V}$ is  the open subset $\{\varphi\in U(X) \ | \
\varphi(x)\in V\}$ of \  $U(X)$.

First we observe that $\Gamma_{x, V}=x^{-1}\Gamma_{e, V}$ and
$$f_*^{-1}(\Gamma_{e, V})\subset f^{-1}(V).$$ Then (\ref{2})
follows from (\ref{1}) and the $X$-equivariance of $f_*$.

Besides, since  $f_*(x)\in \Gamma_{x, V}$  for every $x\in X$, we see that the sets $\Gamma_{x, V}$, $x\in X$, constitute a covering  of $Z$.

\medskip

From now on we restrict ourselves only by the induced actions of
the subgroup $G\subset X$, i.e., we will consider $X$ and $Z$ just as
$G$-spaces.

\medskip

Now, it follows from (\ref{2}) and from the $G$-equivariance of $f_*$ that
 \begin{equation}\label{3}
 f_*^{-1}\big(G(\Gamma_{x, V})\big)\subset UxG, \quad \text{for every} \quad  x\in X.
\end{equation}

Since $f_*:X\to Z$  is $G$-equivariant, it induces a continuous map
$\widetilde{f_*}$ of the $G$-orbit spaces, i.e., we have the following
commutative diagram:

$$\begin{array}{cccccccccc}
X &&\stackrel{f_*}{\longrightarrow} &&  Z \\

\\

\downarrow\lefteqn {p}&&&& \downarrow\lefteqn {q}\\

\\

X/G&& \stackrel{\widetilde{f_*}}{\longrightarrow}&& Z/G
\end{array}
$$
where $p$ and $q$ are the $G$-orbit maps.

It follows from (\ref{3}) that
 \begin{equation}\label{4}
 \widetilde{f_*}^{-1}\big(q(\Gamma_{x, V})\big)\subset  p(Ux), \quad\text{for every}\quad x\in X.
\end{equation}

Thus, the open covering $\{p(Ux) \ | \ x\in X\}$ of the
$G$-orbit space $X/G$ is refined by the open covering $\{
\widetilde{f_*}^{-1}\big(q(\Gamma_{x, V})\big)  \ | \
x\in X\}$.

Since  the metric of $Z$ is $G$-invariant, the orbit space $Z/G$
is pseudometrizable (see Preliminaries). Hence the open covering
$\{q(\Gamma_{x, V}) \ | \ x\in X\}$ of $Z/G$ admits an open
locally finite refinement, say $\{W_i\ | \ i\in \mathcal I\}$ (see
\cite[Theorem~4.4.1 and Remark~4.4.2]{eng:77}). Then, clearly,
$\{p^{-1}\big(\widetilde{f_*}^{-1}\big(W_i)\big) \ | \ i\in
\mathcal I\}$ is an open locally finite refinement of $\{UxG \ | \
x\in X\}$ consisting of $G$-invariant sets. It then follows that
each set $p^{-1}\big(\widetilde{f_*}^{-1}\big(W_i)\big)$ is
contained  in some set $UxG$, $x\in X$, which yields that
$$p^{-1}\big(\widetilde{f_*}^{-1}\big(W_i)\big)=S_iG,$$
where $S_i=p^{-1}\big(\widetilde{f_*}^{-1}\big(W_i)\big)\cap Ux$.

Now, $S_i$,  being a subset of the $G$-small set $Ux$ (see Lemma~\ref{P:right}), is itself  $G$-small. Thus $\{S_iG \mid i\in
\mathcal I\}$ is the desired covering.
\end{proof}

\noindent
{\it Proof of Theorem\ref{T:12}}.
Let  $\{S_iG \ |\ i\in \mathcal I\}$ be the  locally finite open covering
of $X$ from Proposition~\ref{P:31}. Then  the union $S=\bigcup\limits_{i\in \mathcal I}
S_i$ is a $G$-small set (see e.g., \cite[Proposition~1.2(d)]{ab:78}). On the other hand,
$$SG=\Big(\bigcup\limits_{i \in \mathcal I} S_i\Big)G=\bigcup_{i \in \mathcal I} S_iG=X,$$
 yielding that $S$ is a $G$-fundamental subset of $X$. Since the closure of a $G$-small set is $G$-small (see e.g., \cite[Proposition~1.2(b)]{ab:78}), the closure  $\overline{S}$ is the desired closed $G$-fundamental set. \qed

\bigskip

\noindent
{\it Proof of Corollary~\ref{C:3}}.  Since $G$ is almost connected, it has a maximal compact subgroup $K$ (see Preliminaries).
 Since by Corollary~\ref{C:2}, the quotient $X/G$ is paracompact,  due to a result of Abels~\cite[Main Theorem]{ab:74}, $X$ admits a global $K$-slice $S$, and hence, it is $G$-homeomorphic  to the twisted product $G\times_KS$ (see Preliminaries). Since the group $K$ is compact, it then follows that the $K$-orbit map $G\times S\to G\times_KS\cong_GX$ is  open and perfect.  This yields immediately that   the restriction  $p|_S: S\to X/G$  of the $G$-orbit map $p: X\to X/G$ is an open and perfect surjection. Now the result follows. \qed

\bigskip

\noindent
{\it Proof of Corollary~\ref{C:5}}.  By virtue of Corollary~\ref{C:0}, $X$ admits a closed subset  $F\subset X$ such that the restriction   of the quotient projection $p:X\to X/G$ to $F$ is a perfect surjection  $p|_F:F\to X/G$. It then follows from the hypothesis that  $F$ has the property $\mathcal P$. Since a locally compact group is paracompact (even, strongly paracompact~\cite[Theorem~3.1.1]{arhtk}), then again by the hypothesis, the product $G\times F$ has the property $\mathcal P$. Since $F$ is a closed $G$-small set, the  action map $G\times F\to X$ is perfect (see  Abels~\cite[Proposition~1.4]{ab:78}), and since $FG=X$ we see  that the map  $G\times F\to X$ is surjective. Then $X$ has the property  $\mathcal P$ by the hypothesis. \qed

\bigskip
\noindent
{\it Proof of Corollary~\ref{C:6}}.  It is  quite similar to the proof of Corollary~\ref{C:3}.\qed

\bigskip

\noindent
{\it Proof of Theorem~\ref{T:13}}.
Since $G$ is almost connected, it has a maximal compact subgroup, say, $K$ (see Preliminaries). Since, by Corollary~\ref{C:2}, the quotient $X/G$ is paracompact, one can apply the Global Slice Theorem  (see Theorem~\ref{T:GSlice}), according to which $X$ admits a global $K$-slice, say $S$. Then  $X$ is $G$-homeomorphic to the twisted product $G\times_KS$ (see Preliminaries). In turn, due to a result of H.~Abels~ \cite[Theorem~2.1]{ab:74}, the twisted product $G\times_KS$ is  homeomorphic  to the product $G/K\times S$. Thus, $X\cong G/K\times S$.

 Since $G$ is locally compact and paracompact (see \cite[Theorem~3.1.1]{arhtk}) and the  quotient map $G\to G/K$ is open and closed, we infer that $G/K$ is also  locally compact and paracompact. Further, since   $S$ is paracompact, according to a theorem of Morita~\cite{mor:53} one has:
\begin{equation}\label{eq1}
 {\rm dim}\, (G/K\times S)\le  {\rm dim}\, G/K + {\rm dim}\, S.
 \end{equation}

 Since $K$ is compact, according to a result of V.~V.~Filippov~\cite{fil:79} one has the  inequality:
\begin{equation}\label{eq2}
 {\rm dim}\, S\le  {\rm dim}\, S/K + {\rm dim}\, q
 \end{equation}
 where $q: S \to S/K$ is the $K$-orbit projection and
 $${\rm dim}\, q=\sup\,\{{\rm dim}\, g^{-1}(a) \ | \ a\in S/K\}.$$
  Further, since $K$ acts freely on $S$, we see that ${\rm dim}\, q ={\rm dim}\, K.$

Consequently, combining (\ref{eq1}) and (\ref{eq2}) one obtains:
\begin{equation}\label{eq4'}
 {\rm dim}\, (G/K\times S)\le  {\rm dim}\, G/K + {\rm dim}\, K +  {\rm dim}\, S/K.
 \end{equation}

Next, since   $X/G\cong (G\times_KS)/G\cong  S/K$ (see Preliminaries) and
$  {\rm dim}\, G / K + {\rm dim}\, K  = {\rm dim}\, G
$
(see Remark~\ref{R}), it then follows from (\ref{eq4'}) that
$$ {\rm dim}\, (G/K\times S)\le  {\rm dim}\, X/G + {\rm dim}\, G,$$
 as required.
\qed

\medskip

\subsection* {Acknowledgement}
The author  thanks the referee for useful comments.

\medskip

\bibliographystyle{amsplain}

\end{document}